\documentclass[11pt]{article}
\usepackage{epsfig}
\usepackage{amssymb}

 \newtheorem{theorem}{Theorem}[section]
 \newtheorem{lemma}[theorem]{Lemma}
 \newtheorem{corollary}[theorem]{Corollary}
 \newtheorem{remark}[theorem]{Remark}
 
\newtheorem{example}[theorem]{Example}
\newtheorem{examples}[theorem]{Examples}
\newtheorem{definition}[theorem]{Definition}
\newtheorem{conjecture}[theorem]{Conjecture}
\newtheorem{exercise}[theorem]{Exercise} 
\newenvironment{proof}{{\it Proof:\/}}{$\Box$\vskip 0.08in}

\def\Label#1{\label{#1}}

\begin{document}
\thispagestyle{empty}
\ 
\vspace{0.5in}
 \begin{center}
 {\LARGE\bf 3-coloring and other elementary invariants of knots\footnote{
%An extended version of the talks given at Colloquium 
%at Odense University, Dec. 9 1993 and at SBF-Seminar in 
%Goettingen, Feb. 21 1994. Slightly extended version of the preprint
%Institut for Mathematik og Datologi, Odense Universitet,
%Preprints 1994, Nr 26, August 1994, ISSN No. 0903-3920
An extended version of two talks given at 
the Mini-semester on Knot Theory at the Stefan Banach International 
Mathematical Center in Warsaw, July 17- August 18 1995.}}
 
\end{center}
\vspace*{0.5in}
 \begin{center}
                      J\'ozef H.~Przytycki  \footnote{ 
Supported by USAF grant 1-443964-22502 while visiting the 
Mathematics Department, U.C. Berkeley.} 
\end{center}
 
\vspace*{0.5in}
Classical knot theory studies the position of a circle (knot) or of several
circles (link) in $R^3$ or $S^3 =R^3\cup \infty$.
The fundamental problem of classical knot theory is the classification
of links (including knots) up to the natural movement in space which 
is called an ambient isotopy. 
To distinguish knots or links we look for invariants of
links, that is, properties of links which are unchanged under ambient isotopy. 
When we look for invariants of links we have to take into account the following
three criteria:
\begin{enumerate}
\item
[1.] Is our invariant easy to compute?
\item
[2.] Is it easy to distinguish elements in the value set of the invariant?

\item
[3.] Is our invariant good at distinguishing links?
\end{enumerate}

The number of components of a link, $com(L)$, is the simplest invariant.
A more interesting link invariant is given by the linking number, 
defined in 1833 by C.F.Gauss using a certain double integral \cite{Ga}.
H.Brunn noted in 1892 that the linking number has a 
simple combinatorial definition \cite{Br}.
\begin{definition}\label{Definition 0.1}
 Let $D$ be an oriented link diagram. Each crossing has an associated
sign: \ $+1$ for
% \ \ \ \ \ \ 
{\psfig{figure=L+maly.eps}}\ 
and  $-1$ for \ 
%\ \ \ \ \ 
{\psfig{figure=L-maly.eps}}.\ 
The global linking number of $D$,
$lk(D)$, is defined to be half of the sum of the signs of crossings between 
different components of the link diagram. If the diagram has no crossings, 
we put $lk(D)=0$.
\end{definition}
To show that a function defined on diagrams of links 
is an invariant of (global isotopy)
of links, we have to interpret global isotopy in terms of diagrams.
This was done by Reidemeister [Re,1927] and Alexander and Briggs
[A-B,1927].

\begin{theorem} [Reidemeister theorem]
\label{Theorem 0.2}
\ \\
Two link diagrams are ambient isotopic if and only if they are connected by a
finite sequence of Reidemeister moves $R_i^{\pm 1}, i=1,2,3$ (see Fig.0.1)
and isotopy (deformation) of the plane of the diagram. The theorem holds also
for oriented links and diagrams. One then has to take into account
all possible coherent orientations of diagrams involved in the moves.
\end{theorem}

\centerline{\psfig{figure=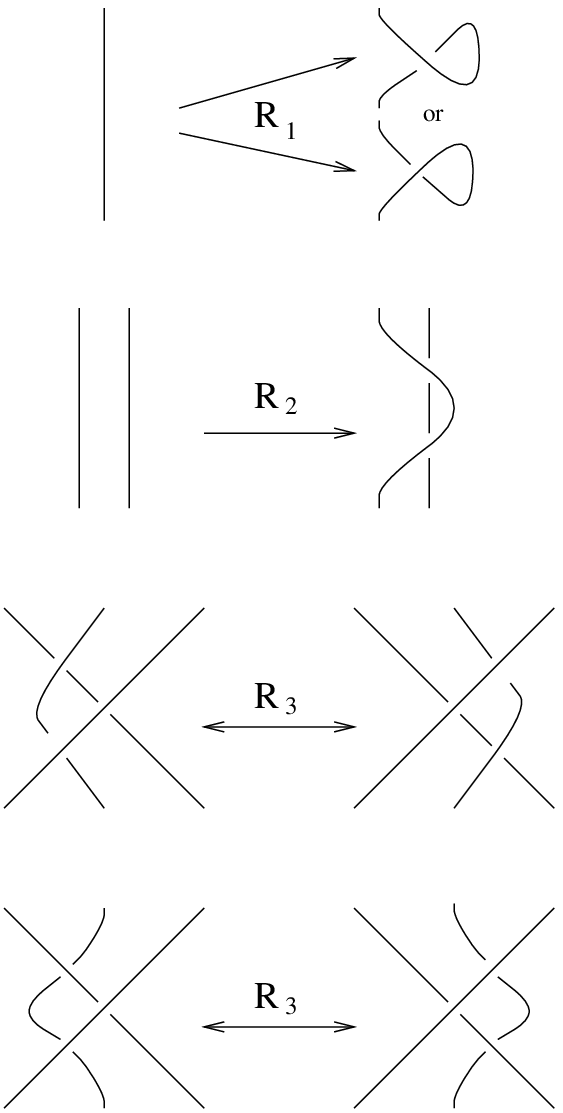}} 
%\vspace*{4.4in}
\begin{center}
Fig. 0.1
\end{center}

\begin{exercise}\label{Exercise 0.3}\ \\
Show that $lk(D)$ is preserved by Reidemeister moves on oriented link diagrams.
Thus $lk$ is an invariant of oriented links.
\end{exercise}
\begin{example}\label{Example 0.4}\ \\
$lk(
{\psfig{figure=T2.eps}}
)=0$, $lk$(
{\psfig{figure=+Hopfmaly.eps}}
)$=1$, 
$lk$(
{\psfig{figure=negHopfmaly.eps}} 
)$=-1$. Therefore the global linking number allows us to distinguish
the trivial link of two components, $T_2$,
%={\psfig{figure=T2.eps}},
the right-handed Hopf link, $2_1$, 
%$2_1=${\psfig{figure=+Hopfmaly.eps}},
and the left-handed Hopf link,  $\bar 2_1$.
%{\psfig{figure=-Hopfmaly.eps}}.
\end{example}

\section{The tricoloring}\label{1}

The tricoloring invariant (or 3-coloring) is the simplest invariant 
which distinguishes between the trefoil knot and the trivial knot.
The idea of tricoloring was introduced by R.Fox 
around{\footnote{Added for e-print: 
I would rather think now that Fox developed the concept 
 around 1956 when he was
explaining Knot Theory to undergraduate students
at Haverford College (``in an attempt to make the subject
accessible to everyone" \cite{C-F}). I am also glad to report other 
articles attempting popularization of Knot Theory to middle and high 
school teachers and students \cite{Cr,Vi,P-6}.}
 1960, [C-F,Chapter VI,Exercises 6-7], \cite{F-2}, and has been extensively
used and popularized by J.Montesinos \cite{Mon}  and L.Kauffman \cite{K}. 
\begin{definition}[{\cite{P-1}}]\label{Definition 1.1}
We say that a link diagram $D$ is tricolored if every arc is colored $r$
(red), $b$ (blue) or $y$ (yellow) ( we consider arcs of the diagram literally,
so that in the undercrossing one arc ends and the second starts; 
compare Fig.1.1), and at any given crossing either all three colors 
appear or only one color appears.
The number of different tricolorings is denoted by $tri(D)$. If a tricoloring
uses only one color we say that it is a trivial tricoloring.
\end{definition}

\centerline{\psfig{figure=trefoil3-col.eps}}
%\vspace*{2in}
\begin{center}
Fig. 1.1. Different colors are marked by lines of different thickness.
\end{center}

\begin{lemma}\label{1.2}
The tricoloring is an (ambient isotopy) link invariant.
\end{lemma}

\begin{proof}

We have to check that $tri(D)$ is preserved under the Reidemeister moves.
The invariance under $R_1$ and $R_2$ is illustrated in Fig.1.2 and the
invariance under $R_3$ is illustrated in Fig.1.3.
\end{proof}
\ \\
\ \\
\centerline{\psfig{figure=R1R23-col.eps}}
%\vspace*{1.4in}
\begin{center}
Fig. 1.2
\end{center}

\centerline{\psfig{figure=R33-col.eps}}
%\vspace*{3in}
\begin{center}
Fig. 1.3
\end{center}

Because the trivial knot has only trivial tricolorings, $tri(T_1)=3$,
and the trefoil knot allows a nontrivial tricoloring (Fig.1.1), 
it follows that the trefoil knot is a nontrivial knot.

\begin {exercise}\label{1.3}
Find the number of tricolorings for the trefoil knot ($3_1$),  the figure 
eight knot ($4_1$) and the square knot ($3_1\# \bar 3_1$, see Fig.1.4). 
Then deduce that these knots are pairwise different. 
\end{exercise}

\begin{lemma}\label{1.4}
$tri(L)$ is always a power of $3$.
\end{lemma}

\begin{proof}
Denote the colors by $0,1$ and $2$ and treat them modulo $3$, that is 
as elements of the group (field) $Z_3$. All colorings of the arcs of 
a diagram using colors $0,1,2$ (not necessarily allowed 3-colorings) 
can be identified with the group $Z_3^r$ where $r$ is the number of arcs 
of the diagram.
The (allowed) 3-colorings can be characterized by the property that at 
each crossing the sum of the colors is equal to zero modulo $3$. 
Thus (allowed) 3-colorings form a subgroup of $Z_3^r$.
\end{proof}
The elementary properties of tricolorings, which we give in Lemma 1.5, 
follow immediately from the connections between tricolorings 
and the Jones and Kauffman polynomials of links.
We will give here an elementary proof of (a)-(c) of Lemma 1.5. 
There is also an elementary proof of (d) (based on the flow-potential 
idea of Jaeger, see \cite{Ja-P}), 
but it is more involved\footnote{Added for e-print: there is another 
elementary proof based on the idea of interpreting tangles as 
Lagrangians in a symplectic space of all 3-colorings of boundary points 
of the tangle \cite{DJP} (see also \cite{P-3,P-4,P-5}).}; compare Lemma 2.2.
\begin{lemma}\label{1.5}
\begin{enumerate}
\item
[(a)]
$tri(L_1)tri(L_2)=3tri(L_1\# L_2)$, where $\#$ denotes 
the connected sum of links
(see Fig. 1.4).
\item
[(b)]
Let $L_+, L_-, L_0$ and $L_{\infty}$ denote four unoriented link diagrams 
as in Fig.1.5. Then, among four numbers $tri(L_+), tri(L_-), tri(L_0)$ and
$tri(L_{\infty})$, three are equal one to another and the fourth is 
equal to them or is $3$ times bigger.\\
in particular:
\item
[(c)] 
$tri(L_+)/tri(L_-)=1$ or $3$, or $1/3$.
\end{enumerate}
Part (b) can be strengthened to show that:
\begin{enumerate}
\item
[(d)] Not all the numbers $tri(L_+), tri(L_-), tri(L_0)$ and $tri(L_{\infty})$
 are equal one to another.
\end{enumerate}
\end{lemma} 
\ \\
\ \\
\centerline{\psfig{figure=connectedsum.eps}}
%\vspace*{2in}
\begin{center}
Fig. 1.4
\end{center}
\ \\
\ \\
 
\centerline{\psfig{figure=L+L-L0Linf.eps}}
%\vspace*{1in}
\begin{center}
Fig. 1.5
\end{center}

\begin{proof}
\begin{enumerate}
\item 
[(a)] An $n$-tangle is a part of a link diagram placed in a 2-disk, with
$2n$ points on the disk boundary ($n$ inputs and $n$ outputs); Fig.1.6.
We show first that for any 3-coloring of a 1-tangle (i.e. a tangle with one
input and one output; see Fig. 1.6(a)), the input arc have the same color
as the output arc.
Namely, let $T$ be our 3-colored tangle and let the 1-tangle $T'$  
be obtained from 
$T$ by adding a trivial component, $C$, below $T$, close to the boundary
of the tangle, so that it cuts $T$ only near the input 
and the output; Fig.1.6(b).
Of course the 3-coloring of $T$ can be extended to a 3-coloring of $T'$ (in
three different ways), because
the tangle $T'$ is ambient isotopic to a tangle obtained from $T$ by adding 
a small trivial
component disjoint from $T$. If we, however, try to color $C$, 
we see immediately that it
is possible iff the input and the output arcs of $T$ have the same color.

Thus if we consider a connected sum $L_1\# L_2$, we see from the above that the arcs joining $L_1$ and $L_2$ have the same color. Therefore the formula,
$tri(L_1)tri(L_2)=3tri(L_1\# L_2)$, follows.
\item
[(b)] Consider a crossing $p$ of the diagram $D$. If we cut out of $D$ a 
neighborhood of $p$, we are left with the 2-tangle, $T_D$ (see Fig.1.6(c)).
The set of 3-colorings of $T_D$, $Tri(T_D)$, forms a $Z_3$ linear space. 
Each of the sets
of 3-colorings of $D_+,D_-,D_0$ and $D_{\infty}$, $Tri(D_+),Tri(D_-),
Tri(D_0)$ and $Tri(D_{\infty})$, respectively form a subspace of 
$Tri(T_D)$. 
Let $x_1,x_2,x_3,x_4$
be generators of $Tri(T_D)$ corresponding to arcs cutting the boundary of 
the tangle; see Fig.1.6(c). 
Then any element of $Tri(T_D)$ satisfies the equality \ $x_1-x_2+x_3-x_4=0$.
To show this, we proceed as in part (a).
Any  element of $Tri(D_+)$ (resp. $Tri(D_-)$,
$Tri(D_0)$ and $Tri(D_{\infty})$) satisfies additionally the equation \
$x_2=x_4$ (resp. $x_1=x_3$, $x_1=x_2$ and $x_1=x_4$). Thus $Tri(D_+)$ (resp.
$Tri(D_-)$, $Tri(D_0)$ and $Tri(D_{\infty})$) is a subspace of $Tri(T_D)$ of
codimension at most one. 
Let $F$ be the subspace of $Tri(T_D)$ given by the equations
$x_1=x_2=x_3=x_4$, that is the space of 3-colorings monochromatic on 
the boundary of
the tangle. $F$ is a subspace of codimension at most one in any 
of the spaces
$Tri(D_+)$, $Tri(D_-)$, $Tri(D_0)$, $Tri(D_{\infty})$. Furthermore
the common part of any two of $Tri(D_+)$, $Tri(D_-),Tri(D_0), Tri(D_{\infty})$
is equal to $F$. To see this we just compare the defining relations 
for these spaces.
Finally notice that
$Tri(D_+)\cup Tri(D_-)\cup Tri(D_0)\cup Tri(D_{\infty})=Tri(T_D)$.
 
We have the following possibilities:
\begin{enumerate}
\item
[(1)] $F$ has codimension 1 in $Tri(T_D)$.\\
Then by the above considerations:\\ 
One of $Tri(D_+),Tri(D_-),Tri(D_0), Tri(D_{\infty})$
is equal to $Tri(T_D)$. The  remaining three spaces are equal to $F$ and (d)
(thus also (b)) of Lemma 1.5 holds.
\item
[(2)] $F=Tri(D_+)=Tri(D_-)=Tri(D_0)= Tri(D_{\infty})=Tri(T_D)$,
\item
[(3)] $F$ has codimension $2$ in $Tri(T_D)$. Then
$3|F|=tri(D_+)=tri(D_-)=tri(D_0)= tri(D_{\infty})=\frac{1}{3}tri(T_D)$
\end{enumerate}
This completes the proof of (b) and (c) of Lemma 1.5. 
To complete (d) of Lemma 1.5 one must exclude cases (2) and (3).
%Case (2) can be excluded relatively easily by showing that $T_D$ has
%a tricoloring nontrivial on the boundary. 
This can be done by showing that tricolorings can be interpreted via 
the so called Goeritz matrix of the link diagram;  
compare Lemma 2.2 and see [J-P]. 

%We consider a link diagram as a four-valent graph (each vertex has 4 edges). 
%Then %consider all colorings of edges of the graph by colors $0,1,2$. Such a
%colorings form a linear space over the field $Z_3$. To have allowed 3-colorings
%we have to put on each crossing two relations.
\end{enumerate}
\end{proof}
\ \\
\ \\
\centerline{\psfig{figure=1-2tangles.eps}}
%\vspace*{1in}
\begin{center}
Fig. 1.6
\end{center}

Part (c) of Lemma 1.5 can be used to approximate the unknotting (Gordian) number
of a knot, $u(K)$; compare \cite{Mur}.

\begin{corollary}\label{Corollary 1.6} \ \\
$u(K)\geq log_3(tri(K))-1$.  In particular
for the square knot:\  $u(3_1\# \bar3_1) =2$.
\end{corollary}

\begin{corollary}\label{Corollary 1.7}
 If $L$ is a link with  $k$-bridge presentation\footnote{
Let $L$ be a link in $R^3$ which meets a plane $E\subset R^3$ in $2k$
points such that the arcs of $L$ contained in each halfspace relative to $E$
possess orthogonal projections onto $E$ which are simple and disjoint.
$(L,E)$ is called a $k$-bridge presentation of $L$; [B-Z].}
then $tri(L)\leq 3^k$.
\end{corollary}

I noticed the connection between tricolorings and the Jones polynomial
when I analyzed the influence of 3-moves on the 3-coloring and 
the Jones polynomial \cite{P-1}.
\begin{definition}\label{Definition 1.8}\ \\
The local change in a link diagram which replaces parallel lines 
by $n$ positive half-twists is called an n-move; see Fig.1.7.
\end{definition}

\begin{lemma}\label{Lemma 1.9}
 Let the diagram $D_{+++}$ be obtained from $D$ by a 3-move
(Fig.1.7(a)). Then:
\begin{enumerate}
\item
[(a)] $tri(D_{+++})=tri(D)$,
\item
[(b)] $V_{D_{+++}}(e^{2\pi i/6})=\pm i^{(com(D_{+++})-com(D))}V_D(e^{2\pi i/6})$, where $V$ is the Jones polynomial, 
\item
[(c)] $F_{D_{+++}}(1,-1)=  F_D(1,-1)$,  where $F$ is the Kauffman polynomial.
\end{enumerate} 
\end{lemma}

\ \\
\centerline{\psfig{figure=n-moves.eps}}
%\vspace*{1in}
\begin{center}
Fig. 1.7
\end{center}
\begin{proof}
We prove (a) and (c) leaving (b) as an exercise.
\begin{enumerate}
\item
[(a)]
The bijection between 3-colorings of $D$ and $D_{+++}$ is illustrated in 
Fig. 1.8.
\ \\ \ \\ \ \\ 
\centerline{\psfig{figure=3-movecol.eps}}
\begin{center}
Fig. 1.8
\end{center} 
\item
[(c)] $F_{D_{+++}}(1,-1)=-F_{D_{+}}(1,-1)-F_{D_{++}}(1,-1)-
F_{D_{\infty}}(1,-1) =
-F_{D_{+}}(1,-1)+ F_D(1,-1)+F_{D_{+}}(1,-1)+F_{D_{\infty}}(1,-1)-
F_{D_{\infty}}(1,-1)= F_D(1,-1)$.
\end{enumerate}
\end{proof}

One can easily check that for a trivial $n$-component
link, $T_n$,
$tri(T_n)=3^n=3V_{T_n}^2(e^{2\pi i/6})=3(-1)^{n-1}F_{T_n}(1,-1)$. Furthermore 
it follows from Lemma 1.9 that as long as a link $L$ can be obtained from a trivial link
by 3-moves we have:
$tri(L)=3|V_L^2(e^{2\pi i/6})|=3|F_L(1,-1)|$.

It may look strange that such a natural problem, whether any link 
can be reduced to an unlink by 3-moves is an open 
problem\footnote{Added for e-print: 
It was showed by M.K.D{\c a}bkowski and the author that Borromean rings 
cannot be reduced by 3-moves to a trivial link. We also found  a smaller 
example -- a closed 3-braid of 20 crossings. The method we develop 
is that of Burnside groups of links which can be interpreted as 
noncommutative version of Fox colorings \cite{D-P-1}.}.
\begin{conjecture} [Montesinos-Nakanishi]
\label{Conjecture 1.10}\ \\
 Any link can be reduced to a trivial link by a sequence of 3-moves.
\end{conjecture}
\begin{remark}\label{Remark 1.11}
Nakanishi first considered the conjecture in 1981. Montesinos
analyzed 3-moves before, in connection with 3-fold dihedral branch
coverings, and asked a related but different question.
Conjecture 1.10 holds for algebraic links (in the Conway sense).
It would be a ``finite" check whether conjecture holds for
links with braid index at most 5 (and bridge index at most 3) as
Coxeter (1957) showed that the quotient of the braid group 
$B_n/<\sigma_1^3>$ is finite for $n\leq 5$.
%[H.S.M. Coxeter, Factor groups of the braid group, Proc. Fourth Canadian
%Math. Congress, Banff, 1957, 95-122.]

According to Nakanishi (1994) the smallest known obstruction to the conjecture
is the 2-parallel of the Borromean rings (notice that it is a 6-string braid),
Fig. 1.9. 
\end{remark}
\ \\
\centerline{\psfig{figure=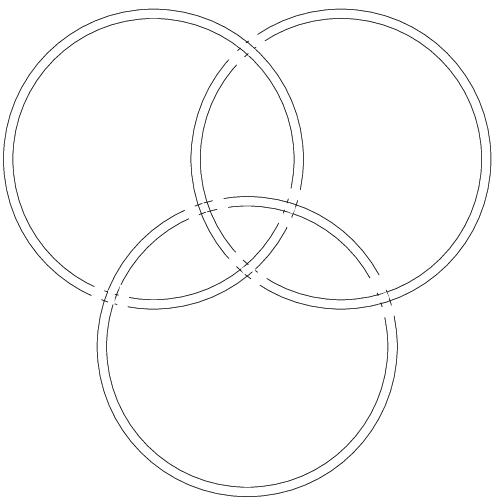}}
\begin{center}
Fig. 1.9 
\end{center}

Lemma 1.5 suggests the following stronger conjecture\footnote{Added 
for e-print: the Conjecture 1.12 does not hold as it is stronger than 
the Montesinos-Nakanishi conjecture which has been 
disproven in \cite{D-P-1}.}.

\begin{conjecture}\label{Conjecture 1.12}\ \\
 Any 2-tangle can be reduced, using 3-moves, to one of the four 2-tangles
of Figure 1.10. We allow additional trivial components in the tangles
of Fig.1.10.
\end{conjecture}
\ \\

\centerline{\psfig{figure=S+S-S0Sinf.eps}}
\begin{center}
Fig. 1.10
\end{center}

Conjecture 1.10 suggests that the formula linking tricoloring with the Jones 
and Kauffman polynomials holds for any link. This is, in fact, the case.

\begin{theorem}\label{Theorem 1.13}\
\begin{enumerate}
\item
[(a)]
$tri(L)=3|V_L^2(e^{2\pi i/6})|$
\item
[(b)]
$tri(L)=3|F_L(1,-1)|$.
\end {enumerate}
\end{theorem}
The proof of (a) in \cite{P-1} uses Fox's interpretation of 3-coloring 
and the connection with the first homology group of the 
branched 2-fold cover of $S^3$ 
branched over the link. Now however we can 
give totally elementary proof based on Lemma 1.5(d).\ \\
\begin{proof}
Because $tri(L)$ is a power of $3$, we can consider the signed version of the
tricoloring defined by:\ \  $tri'(L)=(-1)^{log_3(tri(L))}tri(L)$.\
It follows from Lemma 1.5 (d) that $$tri'(L_+)+tri'(L_-)=
-tri'(L_0)-tri'(L_{\infty}).$$
This is however exactly the recursive formula for the Kauffman polynomial 
$F_L(a,x)$ at $(a,x)=(1,-1)$. 
Comparing the initial data (for the unknot) of $tri'$ and $F(1,-1)$
we get generally that:  \ $-3F_L(1,-1)=tri'(L)=(-1)^{log_3(tri(L))}tri(L)$, 
which proves part (b) of Theorem 1.13. 
Part (a) follows from Lickorish's observation \cite{Li}, that 
$F_L(1,-1)=(-1)^{com(L)}V^2(e^{2\pi i/6})$.
\end{proof}
The value of the Jones polynomial $V_L(e^{2\pi i/6})$ is a slightly more 
delicate invariant than the tricoloring, $tri(L)$, or $F_L(1,-1)$ 
(essentially it is just a ``sign").
P.~Traczyk has given, however, an idea\footnote{Added for 
e-print: see \cite{Tr}.} which allows us to utilize this sign to 
approximate the unknotting number in a better way than in Corollary 1.6.

\begin{theorem}\label{Theorem 1.14}
Let $r(L)=log_3|V^2_L(e^{2\pi i/6})|$ $(=log_3(tri(L))-1)$. 
If a knot $K$ can be trivialized by changing $r(K)$ crossings and  
$V_K(e^{2\pi i/6})={\epsilon}_K (i \sqrt 3)^{r(K)}$, where 
${\epsilon}_K =\pm 1$, \ then the number of negative
crossings, which are changed, is congruent to  $l({\epsilon}_K)$ modulo $2$, 
where $(-1)^{l({\epsilon}_K)}= {\epsilon}_K$.
\end {theorem}
\begin{proof}
Let $t^{1/2}=-e^{2\pi i/12}$, then for any link $L$:\ $V_L(e^{2\pi i/6})=
{\epsilon}_Li^{com(L)-1} (i \sqrt 3)^{r(L)}$. Consider a pair of oriented
links $L_+$ and $L_-$. From the skein relation of the Jones polynomial,
one gets:
$$\frac{1}{2}((1-i\sqrt{3}){\epsilon}_{L_+}i^{com(L_+)-1} (i\sqrt 3)^{r(L_+)}-
(1+i\sqrt{3}){\epsilon}_{L_-}i^{com(L_-)-1} (i\sqrt 3)^{r(L_-)})=$$
$$-i{\epsilon}_{L_0}i^{com(L_0)-1} (i\sqrt 3)^{r(L_0)}$$
One can see immediately that the above equation cannot hold for
$|r(L_+)-r(L_-)|\geq 2$. For $r(L_+)-r(L_-)=1$ it simplifies to:\\
$\frac{1}{2}((3+i\sqrt{3}){\epsilon}_{L_+}-(1+
i\sqrt{3}){\epsilon}_{L_-})=
(-1)^{(\frac{com(L_0)-com(L_+)-1}{2})}{\epsilon}_{L_0}
(i\sqrt{3})^{(r(L_0)-r(L_-))}$.
This equation holds iff
${\epsilon}_{L_+}={\epsilon}_{L_-}$.
Similarly, for $r(L_+)-r(L_-)=-1$ one gets ${\epsilon}_{L_+}=-{\epsilon}_{L_-}$.
This completes the proof of Theorem 1.14, because for the trivial knot, $T_1$,
${\epsilon}_{T_1}=1$.
\end{proof}

\begin{examples}\label{Example 1.15}
\begin{enumerate}
\item
[(a)] Let $K= 3_1\#\bar{3}_1$, then $u(K)=2$. Furthermore if $K$ is trivialized
using two crossing changes, then one positive and one negative crossing, 
have to be changed. Namely $V_K(e^{2\pi i/6})=3=-(i\sqrt 3)^2$ and Theorem 1.14 
can be used.
\item
[(b)] The unknotting number of the knot $7_7$ is $1$; see Fig.1.11.
However this knot cannot be trivialized by changing a positive crossing.
Namely  $V_{7_7}(e^{2\pi i/6})=-i\sqrt{3}$. Notice that the signature of
$7_7$ is equal to $0$ and the Tait number of the minimal diagram is equal to
$+1$.
\end{enumerate}
\end{examples}
\ \\ 
\centerline{\psfig{figure=7-7trivialization.eps}}
%\vspace*{2in}
\begin{center}
Fig. 1.11
\end{center}

\section{n-coloring}\label{2}

Tricoloring of links can be generalized, after Fox, 
[F-1;Chapter 6], [C-F;Chapter VIII,Exercises 8-10], \cite{F-2}, 
to n-coloring of links as follows:
\begin{definition}\label{Definition 2.1}
We say that a link diagram $D$ is n-colored if every
arc is colored by one of the numbers $0,1,...,n-1$ in such a way that 
at each crossing the sum of the colors of the undercrossings is equal 
to twice the color
of the overcrossing modulo $n$. 
\end{definition}
The following properties of n-colorings, can be proved in a similar way
as the tricoloring properties. However, an elementary proof of the part (g)
is more involved and requires an interpretation of n-colorings using the
Goeritz matrix \cite{Ja-P}.
\begin{lemma}\label{2.2}
\begin{enumerate}
\item
[(a)] Reidemeister moves preserve the number of n-colorings, $col_n(D)$, 
      thus it is a link invariant,
\item
[(b)] if $D$ and $D'$ are related by a finite sequence of $n$-moves, then 
$col_n(D)=col_n(D')$,
\item
[(c)] $n$-colorings form an abelian group, $Col_n(D)$,
\item
[(d)] if $n$ is a prime number, then $col_n(D)$ is a power of $n$ and for a
link with $b$ bridges: \ \ $b\geq log_n(col_n(L))$,
\item
[(e)] $col_n(L_1) col_n(L_2)=n(col_n(L_1\# L_2))$,
\item
[(f)] if $n$ is a prime odd number then among the four numbers 
$col_n(L_+), col_n(L_-),$ $col_n(L_0)$ and $col_n(L_{\infty}),$ 
three are equal one to another and the fourth is 
either equal to them or $n$ times bigger,
\end{enumerate}
More generally: If $L_0, L_1, ..., L_{n-1},L_{\infty}$ are $n+1$ diagrams
generalizing the four diagrams from (f); see Fig.2.1 then:
\begin {enumerate}
\item
[(g)] if $n$ is a prime number then among the $n+1$ numbers 
$col_n(L_0), col_n(L_1),...,col_n(L_{n-1})$ and $col_n(L_{\infty})$ 
$n$ are equal one to another and the $(n+1)$'th is $n$ times bigger,
\item
[(h)]
if $n$ is a prime number, then $u(K)\geq log_n(col_n(K))-1$.
\end{enumerate}  
\end{lemma}
\ \\

\centerline{\psfig{figure=L0L1L2L3Linf.eps}}
%\vspace*{2in}
\begin{center}
Fig. 2.1
\end{center}

\begin{corollary}\label{Corollary 2.3}
% When $col_n(L)$ is interpreted using homology of double cover then
% we use Goeritz matrix modulo $n$.
\begin{enumerate}
\item
[(i)] For the figure eight knot, $4_1$, one has $col_5(4_1)=25$, so the figure
eight knot is a nontrivial knot; compare Fig.2.2.
\item
[(ii)] $u(4_1\# 4_1)=2$.
\end{enumerate}
\end{corollary}

\centerline{\psfig{figure=4-1fivecol.eps}}
%\vspace*{2in}
\begin{center}
Fig. 2.2
\end{center}

By Lemma 2.2(b), any 5-move preserves the number of 5-colorings. On the other
hand, Corollary 2.3 suggests that the move of Fig.2.3 also 
preserves $col_5(L)$.
\begin{lemma}\label{Lemma 2.4}\ If two links are related by a sequence of moves
as in Fig.2.3 (allowing the mirror image of Fig.2.3), then they have the 
same number of 5-colorings.
\end{lemma}

\ \\

\centerline{\psfig{figure=claspmove.eps}}
%\vspace*{2in}
\begin{center}
Fig. 2.3
\end{center}

\begin{proof} It suffices to notice that for $x$ and $y$ of Fig.2.3:\\
$3x-2y \equiv \frac{x+y}{2}$ mod $5$ and $2x-y \equiv \frac{3y-x}{2}$ 
mod $5$.
\end{proof}

It was noticed in \cite{H-U}, that the moves of Fig.2.3 are more general
than the 5-moves.
\begin{lemma}[{\cite{H-U}}] \label{Lemma 2.5} \
A 5-move is a combination of moves of Fig.2.3 (and isotopy).
\end{lemma}
\begin{proof}
This is illustrated in Fig.2.4.
\end{proof}
\ \\

\centerline{\psfig{figure=clasp-5moves.eps}}
%\vspace*{2in}
\begin{center}
Fig. 2.4
\end{center}

It has been noticed in \cite{P-2} that there are links which cannot be changed
to trivial links using 5-moves. In particular, the figure eight 
 knot cannot be reduced to a trivial link by
5-moves (this is an easy application of the Jones polynomial evaluated
at $t=e^{\pi i/5}$).
 It is however an open problem whether
any link can be changed to a trivial link by moves of the type shown in
Fig.2.3 \cite{Nak} (it holds for links up to 7 crossings\footnote{Added 
for e-print; It holds for links up to 8 crossings but the knot $9_{49}$ 
cannot be reduced to a trivial link by these moves \cite{D-P-2}.}, 
in particular for the Borromean rings). More generally, it holds
for algebraic knots (in the sense of Conway).\\
The immediate generalization of the move of Fig.2.3 is a move which
changes $p$ horizontal half twists into $q$ vertical half twists.
Let us call such a move (and its mirror image) a $[p,q]$-move; see
Fig. 2.5. 
\ \\
\centerline{\psfig{figure=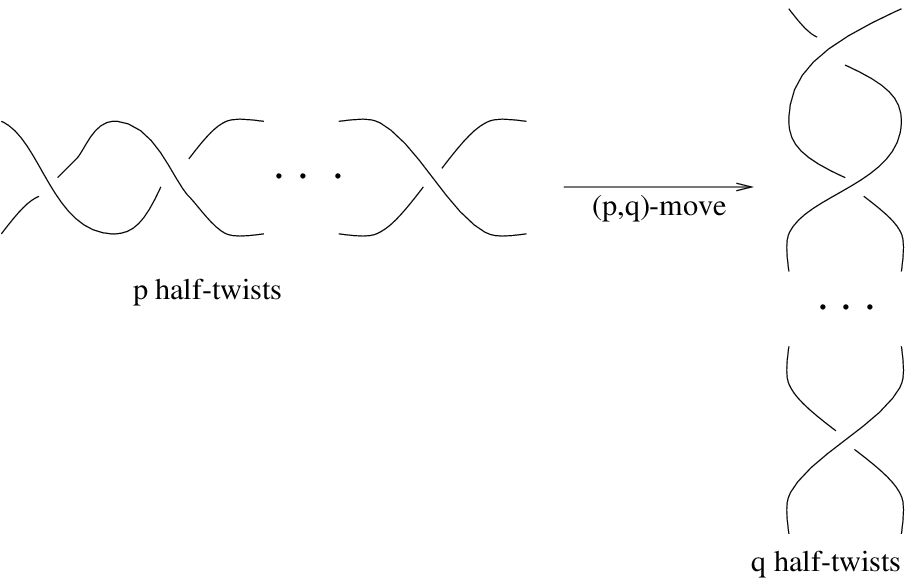}}
\begin{center}
Fig. 2.5
\end{center}\
\ \\

\begin{exercise}
\begin{enumerate}
\item
[(a)]
Show that a $[p,q]$-move preserves the number of $(pq+1)$-colorings.
\item [(b)]
Show that a $(2p+1)$-move is a combination of a $[p,2]$-move and 
a $[2,p]$-move. 
\end{enumerate}
It is not always true that a $(pq+1)$-move is a composition  of 
$[p,q]$-moves (and their inverses).
\begin{enumerate}
\item
[(c)]
Use linking numbers to show that if $p,q$ are odd numbers 
and $p+q$ is not a divisor of $pq+1$ (e.g. $p=q\geq 3$), then a 
 $(pq+1)$-move is not a composition of $[p,q]$-moves
(and their inverses).\ Show in particular that the torus link of type 
$(10,2)$ is not $[3,3]$ equivalent to the trivial link of 2 components
(that is the torus link cannot be obtained from the trivial link by
the sequence of $[3,3]$-moves).
% a $[p,q]$-move preserves linking numbers modulo $\frac{p+q}{2}$, 
%but $pq+1$-move may change it.  
\item
[(d)] Use the Kauffman polynomial to show that the $17$-move is not a
composition of $[4,4]$-moves.
%Show that a $(p(p+1)+1)$ is not a composition of $[p,p+1]$-moves
%(and their inverses), for $p>2$.
\end{enumerate}
\end{exercise}
{\it Hint to (d)}. Analyze how the Brandt-Lickorish-Millett polynomial,
$Q_L(x)=F(1,x)$  changes under $[p,q]$-moves and $(pq+1)$-moves; 
compare \cite{P-2}. Figure 2.6 illustrates the fact that a $[4,4]$
move preserves, up to the factor $-1$, the $Q_L(x)$ polynomial 
modulo $x^4+x^3-2x^2+1$.  A $17$-move can change the polynomial $|Q_L(x)|$.

\ \\
\centerline{\psfig{figure=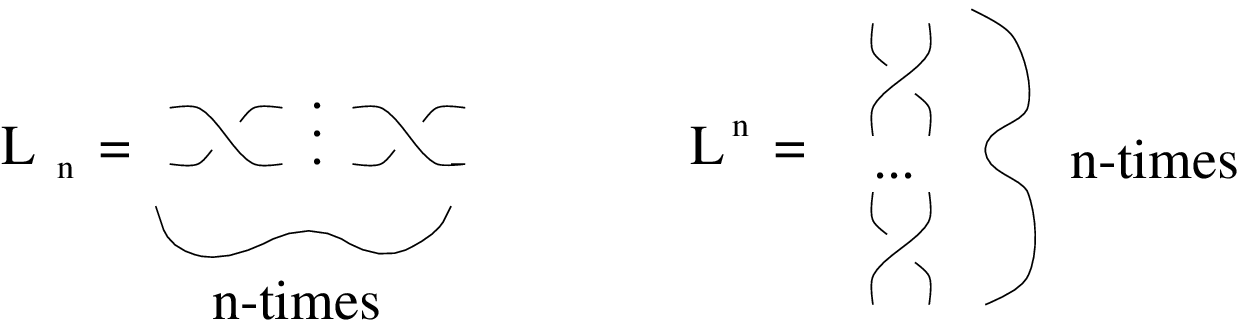}}
\begin{center}
Fig. 2.6
\end{center}\ \ \\

 Consider the Brandt-Lickorish-Millett polynomial of
$L_4(x)$ and $L^4(x)$. From the relation $Q_{L_n}=xQ_{L_{n-1}}-
Q_{L_{n-2}} +x Q_{L_{\infty}}$ one obtains\footnote{It can be easily 
checked by induction, 
see \cite{P-2}, that generally 
$Q_{L_n}= U_{n-1}Q_{L_1} - U_{n-2}Q_{L_0} + \frac{x}{x-2}
(U_{n-1}-U_{n-2}-1)Q_{L_{\infty}}$ where $U_i(x)$ is the Chebyshev polynomial
of the second type defined by: $U_0(x)=1$, $U_1(x)=x$, 
$U_i(x) = xU_{i-1}(x)-U_{i-2}(x)$.}\\ 
$Q_{L_4}= (x^3-2x)Q_{L_1}- (x^2-1)Q_{L_0} + (x^3 + x^2)Q_{L_{\infty}}$,\ 
and\\
 $Q_{L^4}=(x^4-3x^2+1)Q_{L^0} - (x^3-2x)Q_{L^{-1}} + 
(x^4+x^3-x^2)Q_{L^{\infty}}=\\ (x^4-3x^2+1)Q_{L_{\infty}} - (x^3-2x)Q_{L_1} +
(x^4+x^3-x^2)Q_{L_0}$. Therefore for $x$ such that $x^4+x^3-2x^2+1=0$,
$Q_{L_4}=-Q_{L^4}$. On the other hand consider the torus link $L_{17,2}$,
which can be reduced to the trivial link of 2 components by a $17$-move. 
If a $17$-move was the combination of
$[4,4]$-moves then $Q_{L_{17,2}}(x)=\pm \frac {2-x}{x}$ 
for $x^4+x^3-2x^2+1=0$. One can check however that $\frac{x}{2}
(Q_{L_{17,2}}(x)+ \frac {2-x}{x})$ is an irreducible polynomial of degree 
equal to $17$ and $Q_{L_{17,2}}(x)- \frac {2-x}{x}=2\frac{x-1}{x}
(1-4x-10x^2+10x^3+15x^4-6x^5-7x^6+x^7+x^8)^2$, which for $x=q+q^{-1}$
gives $2\frac{q-1+q^{-1}}{q+q^{-1}}q^{-16}(\frac{q^{17}-1}{q-1})^2$. Thus
$L_{17,2}$ is not $[4,4]$ equivalent to the trivial link of 2 components.
\ \\
%\begin{center} Fig. 2.6 \end{center}\ 

It is a nice exercise in linear algebra to show that the number of
 n-colorings is preserved by certain generalizations of n-moves.\\
Let $t_{\Delta ,k}$ denote the righthanded half-twist performed on k strings;
see Fig.2.7.\\
\ \\
\centerline{\psfig{figure=tDeltak.eps}}
\begin{center}
Fig. 2.7
\end{center}\
\begin{lemma}\label{Lemma 2.7}
\begin{enumerate}
\item [(a)]
$t^4_{\Delta ,k}$ preserves $col_n(D)$, for odd $k$ and any $n$.
\item
[(b)] $t^{2n}_{\Delta ,k}$ preserves $col_{2n}(D)$, for an even $k$.
\item [(c)]
Lemma 2.2(b) is stronger than (b) for $k=2$, and can be written as:\\
$t^n_{\Delta ,2}$ preserves $col_n(D)$.
\end{enumerate}
\end{lemma}

\section{Coloring and algebraic topology}\label{3}

It is useful to look at Lemma 2.2 from the point of view of algebraic
topology.
\begin{definition}\Label{Definition 3.1}
Consider the abelian group of all colorings of arcs of a diagram using 
integers as colors. In other words consider the free abelian group spanned 
by all arcs of the diagram. Let each crossing give the relation: 
the sum of the colors of the undercrossings is equal to twice the color 
of the overcrossing. Let $H_D$ denote the
described group.
\end{definition}
\begin{lemma}\label{Lemma 3.2}\ \
\begin{enumerate}
\item
[(a)] $H_D$ is preserved by Reidemeister moves, therefore it is
a link invariant, $H_L$,
\item
[(b)] $H_D$ reduced modulo $n$ (i.e. $H_D\otimes Z_n$) is the group of 
$n$-colorings of $D$.
\end{enumerate}  
\end{lemma}
\begin{theorem}\label{Theorem 3.3}
$H_L$ is the direct sum of the first homology group of the cyclic 
branched double cover of $S^3$ with branching set $L$ and 
the infinite cyclic group. That is:\ \ $H_L=H_1((M_L)^{(2)},Z)\oplus Z$.
\end{theorem}
Before we offer two proofs of the theorem we can carry our combinatorial 
construction one step further.
\begin{definition}[{\cite{F-R}}]\label{Definition 3.4}\ \\
$G_D$ is the group associated to the diagram $D$ as follows: \ \
generators of $G_D$ correspond to arcs of the diagram. Any crossing $v_s$ yields
the relation $r_s=y_iy_j^{-1}y_iy_k^{-1}$ where $y_i$ corresponds to the overcrossing
and $y_j,y_k$ correspond to the undercrossings at $v_s$.
\end{definition}
The group $G_D$ was introduced by R.Fenn and C.Rourke \cite{F-R} as an example 
of a rack's functor. They call it the associated core group of a link; 
compare with the core group of Joyce \cite{Joy}, an example of an involutory
quandle. Joyce refers to the 1958 book of Bruck \cite{Bruc} as the source of 
the idea; compare paragraphs 1, 2 and 19 of \cite{Joy}.
The topological interpretation of $G_D$ was given by M.Wada \cite{Wa}; 
see Theorem 3.6.
\begin{lemma}\label{Lemma 3.5}
\begin{enumerate}
\item 
[(a)] $G_D$ is preserved by Reidemeister moves, thus it is 
a link invariant, $G_L$, 
\item 
[(b)] The abelianization of $G_D$ yields $H_D$.
\end{enumerate}   
\end{lemma} 
Lemma 3.5(b) and Theorem 3.3  suggest that $G_D$ may be related to the
fundamental group of the branched 2-fold cover over $D$. This is 
in fact the case.\\
\ \\
\ \\
\begin{theorem}[{\cite{Wa}}]\label{Theorem 3.6}\ \\ 
$G_L$ is the free product of the fundamental group of the cyclic 
branched double
cover of $S^3$ with branching set $L$ and the infinite cyclic group.
That is:\ \ $G_L={\pi}_1((M_L)^{(2)})\ast Z$. 
\end{theorem} 
We will give later an elementary proof of the Wada theorem.

To prove Theorems 3.3 and 3.6, we need a combinatorial definition of another, 
well known, group.
\begin{definition}\label{Definition 3.7}
${\Pi}_D$ is the group associated to an oriented link diagram $D$ as follows:\ \
generators of the group, $x_1,...,x_n$, correspond to arcs of the diagram; any
crossing $v_s$ yields the relation $r_s=x_i^{-1}x_j^{-1}x_ix_k$  or
$r_s= x_ix_jx_i^{-1}x_k^{-1}$,
where $x_i$ corresponds to the overcrossing
and $x_j,x_k$ correspond to the undercrossings at $v_s$ and the first relation
comes from a positive crossing, Fig.3.1(a), and the second comes from a
negative crossing, Fig.3.1(b).
\end{definition}
\ \\
\ \\
\ \\
\centerline{\psfig{figure=Wirtinger.eps}}
%\vspace*{1.4in}
\begin{center}
Fig. 3.1 
\end{center}

\begin{lemma}\Label{Lemma 3.8}
\begin{enumerate}
\item
[(a)]
${\Pi}_D$ is preserved by Reidemeister moves, thus it is
a link invariant,
\item
[(b)] Abelianization of the group ${\Pi}_D$ is a free group of $com(D)$ 
generators.
\item
[(c)] ${\Pi}_D$ does not depend on the orientation of the link. In particular 
if we change the orientation of a component, say $D_1$, of $D$ to get 
the diagram $D'$, then the isomorphism of the group ${\Pi}_D$ generated 
by $(x_1,...,x_n)$ onto the group
${\Pi}_{D'}$ generated by $(x'_1,...,x'_n)$ is given 
by sending $x_i$ to $x_i'$ or ${x'}_i^{-1}$ depending on whether the arc 
of $x_i$ preserves or changes orientation when going from $D$ to $D'$.
\end{enumerate}
\end{lemma}
The group ${\Pi}_D$ and its presentation, which we described, was introduced by
W.~Wirtinger at his lecture delivered at a meeting of the German 
Mathematical Society in 1905 \cite{Wi}.

\begin{theorem}[Wirtinger]\label{Theorem 3.9}\ \\
${\Pi}_D$ is the fundamental group of the complement of the link; i.e.
${\Pi}_D={\pi}_1(S^3-D)$.
\end{theorem}
We will not use this theorem, except for further algebraic-topological 
interpretations. For the proof see \cite{C-F}, \cite{Rol}, or \cite{B-Z}.

With the group ${\Pi}_D$ defined, we can give Fox's interpretation 
of $n$-colorings.

Let $D_n$ denote the dihedral group, i.e. the group of isometries of a regular
$n$-gon. $D_n$ has a presentation: \
 $D_n=\{\alpha ,s:\ {\alpha}^n=1,\ s^2=1,\ s\alpha s={\alpha}^{-1}\}$. 
The rotations, $\{{\alpha}^k\}$,  form a cyclic subgroup, a $Z_n$. 
Reflections can be written as: \ $s_k=s{\alpha}^k$.
\begin{lemma}\label{Lemma 3.10}
$n$-colorings of $D$ are in bijection with homomorphisms 
from ${\Pi}_D$ to $D_n$, 
which send $x_i$ to reflections. 
Namely, for an $n$-coloring $c$, the homomorphism
${\phi}_c:{\Pi}_D\to D_n$ is given by: \ ${\phi}_c(x_i)=s_k$, where
$k$ is the color of the arc which correspond to $x_i$.
\end{lemma}
To prove Theorem 3.3 we need still more preparation.

Let $\nu :{\Pi}_D\to Z_2$ be the modulo $2$ evaluation map, that is it sends
words of even length (in the generators $x_i^{\pm 1}$) to $0$, and words of odd
length to one. $\nu$ is well defined because the relations of the group
${\Pi}_D$ have even length. Denote by ${\Pi}_D^{(2)}$ the kernel, $ker(\nu )$,
of the epimorphism $\nu$. This is a subgroup of index $2$ in ${\Pi}_D$.
From the point of view of algebraic topology it is the fundamental group of
the 2-fold cyclic covering of $S^3-D$. That is
${\Pi}_D^{(2)}={\pi}_1(({S^3-D})^{(2)})$. The abelianization of ${\Pi}_D^{(2)}$
is the first homology group of $({S^3-D})^{(2)}$.
\begin{lemma}\label{3.11}
\begin{enumerate}
\item
[(a)] For any $n$-coloring, $c$, one has ${\phi}^{-1}_c(Z_n)={\Pi}_D^{(2)}$.
\item
[(b)] Any homomorphism ${\phi}_c:{\Pi}_D\to D_n$ lifts uniquely to 
a homomorphism \\
${\phi}^{(2)}_c:{\Pi}_D^{(2)} \to Z_n$. 
In particular ${\phi}^{(2)}_c(x_i^2)=0$.
\item
[(c)] For any $n$-coloring, $c$, the following diagram is commutative: \ \\
$\begin{array}{ccc}
{\Pi}_D^{(2)} & \stackrel{{\phi}_c^{(2)}}{\longrightarrow} & Z_n  \\
\downarrow    &                 & \downarrow \\
{\Pi}_D       &  \stackrel{{\phi}_c}{\longrightarrow} & D_n  \\
\downarrow    &                 & \downarrow \\
Z_2           &  =              & Z_2
\end{array}$
\item
[(d)] Any homomorphism ${\phi}^{(2)}:{\Pi}_D^{(2)} \to Z_n$, such that
 ${\phi}^{(2)}(x_i^2)=0$ is a lift of exactly $n$ homomorphisms
${\phi}_c:{\Pi}_D\to D_n$. 
\end{enumerate}
\end{lemma}
\begin{proof}
\begin{enumerate}
\item
[(a)] This is the case because rotations in $D_n$ (i.e. $Z_n$) 
are compositions of an even number of reflections.
\item 
[(b)] This follows from (a) and reflects the fact that ${\phi}_c$ sends 
words of even length to words of even length.
\item
[(c)] This summarizes (a) and (b).
\item
[(d)] Fix $x_i$. To show (d) it suffices to show that for any $s_j$, 
there is exactly one $c$ such that 
${\phi}_c:{\Pi}_D\to D_n$ lifts to ${\phi}^{(2)}$ and ${\phi}_c(x_i)= s_j$. 
Namely we define ${\phi}_c(w)={\phi}^{(2)}(w)$ if $w$ has an even length, and
${\phi}_c(w)={\phi}^{(2)}(wx_i^{-1})s_j$ if $w$ has an odd length. 
In particular ${\phi}_c(x_k)={\phi}^{(2)}(x_kx_i^{-1})s_j$. 
We have to check that ${\phi}_c$ is
a homomorphism. Consider ${\phi}_c(w_1){\phi}_c(w_2)$. We have to check four
cases, however for $w_1$ of even length the checking is immediate (e.g. if
$w_2$ has odd length then ${\phi}_c(w_1){\phi}_c(w_2)={\phi}^{(2)}(w_1)
{\phi}^{(2)}(w_2x_i^{-1})s_j={\phi}^{(2)}(w_1w_2x_i^{-1})s_j=
{\phi}_c(w_1w_2)$). For the other cases we need to check first that 
${\phi}^{(2)}(ww)=0$ for any $w$ of odd length. We will show 
the slightly stronger fact
that $ww$ lies in the commutator subgroup of the quotient group 
$\bar{\Pi}_D^{(2)}={\Pi}_D^{(2)}/(x_i^2)$. 
Namely, let $w=x_{i_1}x_{i_2}...x_{i_{2m+1}}$,
then $ww=x_{i_1}x_{i_2}...x_{i_{2m+1}}x_{i_1}x_{i_2}...x_{i_{2m+1}}=$ 
$(x_{i_1}x_1)(x_1x_{i_2})...
(x_{i_{2m+1}}x_1)(x_1x_{i_1})...$\
$(x_{i_{2m}}x_1)
( x_1x_{i_{2m+1}})=$
$(x_{i_1}x_1)(x_{i_2}x_1)^{-1}
(x_{i_3}x_1)(x_{i_4}x_1)^{-1}...$\ $(x_{i_{2m+1}}x_1)(x_{i_1}x_1)^{-1}...$ 
$(x_{i_{2m+1}}x_1)^{-1} $. \ \ Now we can check that ${\phi}_c(w_1){\phi}_c(w_2)=
{\phi}_c(w_1w_2))$ for $w_1$ of odd length.
\begin{enumerate}
\item
[(i)] If $w_2$ is of odd length then:\\
${\phi}_c(w_1){\phi}_c(w_2)=
{\phi}^{(2)}(w_1x_i^{-1})s_j{\phi}^{(2)}(w_2x_i^{-1})s_j=$
${\phi}^{(2)}(w_1x_i^{-1})({\phi}^{(2)}(w_2x_i^{-1}))^{-1}=
{\phi}^{(2)}(w_1x_i^{-1}) {\phi}^{(2)}(x_iw_2^{-1})=$
${\phi}^{(2)}(w_1w_2^{-1})={\phi}^{(2)}(w_1w_2){\phi}^{(2)}(w_2^{-2})=
{\phi}_c(w_1w_2))$
\item
[(ii)] If $w_2$ is of even length, then using (i) we get:\\
${\phi}_c(w_1){\phi}_c(w_2)={\phi}_c(w_1){\phi}^{(2)}(w_2x_i^{-1}x_i)=
{\phi}_c(w_1){\phi}_c(w_2x_i^{-1}){\phi}_c(x_i)=$
${\phi}^{(2)}(w_1w_2x_1^{-1})s_j=
{\phi}_c(w_1w_2)$.
\end{enumerate}
\end{enumerate}
\end{proof}
The quotient group $\bar{\Pi}_D^{(2)}={\Pi}_D^{(2)}/(x_i^2)$
can be interpreted as the fundamental group of the cyclic 
branched double cover of $S^3$ with branching set $D$; that is 
$\bar{\Pi}_D^{(2)}={\pi}_1((M_D)^{(2)})$. 
This interpretation follows from the fact that the elements $x_i^2$ 
correspond to meridians of boundary components of the unbranched
double cover of $S^3-D$ and that these meridians are ``killed" in the branched
cover.  The homomorphism ${\phi}^{(2)}$, from Lemma 3.11(d), factors through 
 $\bar{\Pi}_D^{(2)}$, and because $Z_n$ is abelian, it factors through the 
abelianization of $\bar{\Pi}_D^{(2)}$. This abelianization can be interpreted
as the first homology group of the cyclic branched double
cover of $S^3$ with branching set $D$. 
We denote this group by $H_1=H_1((M_D)^{(2)},Z)$.
Therefore we have the following commutative diagram:\ \\
$\begin{array}{ccl}
\bar{\Pi}_D^{(2)} & \longleftarrow  & {\Pi}_D^{(2)} \\
\downarrow        &                 & \downarrow {\phi}^{(2)} \\
H_1               & \longrightarrow & Z_n
\end{array}$ \ \\
We have also a bijection between homomorphisms $H_1 \to Z_n$, and homomorphisms
${\phi}^{(2)}:{\Pi}_D^{(2)} \to Z_n$, which satisfy condition (d) of Lemma 3.11.
Thus Lemma 3.11 (d) leads to:
\begin{corollary}\label{3.12}
To any homomorphism $H_1\to Z_n$, there is uniquely associated $n$ different 
$n$-colorings. In particular the trivial homomorphism corresponds to 
$n$ trivial $n$-colorings. Therefore\ \
$H_1\otimes Z_n \oplus Z_n =
Hom(H_1\oplus Z_n,Z_n)$ has the same number of elements as $H_D\otimes Z_n$.
\end{corollary}
Because Corollary 3.12  holds for any $n$, we have $H_D=H_1\oplus Z$ 
and the proof of Theorem 3.3 is completed.

\begin{lemma}\label{Lemma 3.13}
Let $F=\{x_1,...,x_n: \}$ be the free group on $n$ generators. 
Let $F^{(2)}$ be its subgroup generated by the words of even length and 
let $\bar{F}^{(2)}$ be the quotient group $F^{(2)}/(x_i^2)$. Then:
\begin{enumerate}
\item
[(a)]
$F^{(2)}$ is a free group on $2n-1$ generators $x_1x_k,x_2x_k,
...,x_{k-1}x_k,x_{k+1}x_k,$ $...,x_nx_k,x_1^2,x_2^2,...,x_n^2$, where $x_k$ is
any fixed generator of $F$.
\item
[(b)]
$\bar F^{(2)}$ is a free group on $n-1$ generators 
$x_1x_k,x_2x_k, ...,x_{k-1}x_k,x_{k+1}x_k,...,x_nx_k$.
\end{enumerate}
\end{lemma}
The above lemma is the starting point of our proof, given below,
 of the Wada theorem (3.6).\\
\begin{proof}
Let $D'$ denote the diagram obtained from the diagram $D$ by adding one
trivial component.
\begin{enumerate} 
\item 
[Step 1]
${\bar \Pi}_{D'}^{(2)}={\bar \Pi}_D^{(2)}\ast Z$.\ \\
\begin{proof}
Denote by $x_1,...,x_n$ generators corresponding to arcs of the diagram $D$, and
by $x_{n+1}$ the generator corresponding to the additional component of $D'$.
${\bar \Pi}_{D'}^{(2)}$ is generated by $x_2x_1, x_3x_1,...,x_nx_1,x_{n+1}x_1$.
Relations of the group are associated to crossings of the diagram $D'$ (so $D$).
The general relation is of the form  
$x_i^{-1}x_jx_ix_k^{-1}$ or $x_ix_jx_i^{-1}x_k^{-1}$ ($i,j,k\leq n$); 
Fig.3.1. Both relations lie in $F^{(2)}$ and in  $\bar F^{(2)}$
they are both conjugated to $(x_ix_1)(x_jx_1)^{-1}(x_ix_1)(x_kx_1)^{-1}$. 
${\bar \Pi}_{D'}^{(2)}$ is a quotient of $\bar F^{(2)}$, 
by these relations (compare Lemma 3.14). No relation
uses the generator $(x_{n+1}x_1).$ Therefore ${\bar \Pi}_{D'}^{(2)}={\bar \Pi}_D^{(2)}\ast Z$.
\end{proof}
\item
[Step 2.] Consider generators $x_1x_{n+1}, x_2x_{n+1},...x_nx_{n+1}$ of 
${\bar \Pi}_{D'}^{(2)}$. We can associate them to arcs of $D$ as follows:\
$x_ix_{n+1}$ corresponds to the arc of $D$ which before was associated to $x_i$. 
No generator corresponds to the additional arc of $D'$. Relations associated
to crossings of $D$ can be found as in Step 1 to be: \ \
$(x_ix_{n+1})(x_jx_{n+1})^{-1}(x_ix_{n+1})(x_kx_{n+1})^{-1}$, where
$i,j,k\leq n$. If we put $y_s=x_sx_{n+1}$ for $s\leq n$ then the relations reduce

to $y_iy_j^{-1}y_iy_k^{-1}$. We get exactly the presentation of the group
$G_D$ from Definition 3.4. Therefore $G_D={\bar \Pi}_{D'}^{(2)}=
{\bar \Pi}_D^{(2)}\ast Z$. The proof of the Wada theorem is completed. 
\end{enumerate}
\end{proof}
One can describe the group 
${\Pi}_D^{(2)}$ similarly as the group ${\bar \Pi}_D^{(2)}$.
A more challenging exercise is to find the Wirtinger type presentation of the
fundamental group of the general $k$ fold cyclic branched cover of $S^3$ with 
branching set $D$, i.e. ${\Pi}_D^{(k)}={\pi}_1(M^{(k)})$.
%In the following lemma we will describe the result for $k=3$.
We describe the result below (compare [B-Z;Ch.4]\footnote{Added for e-print: 
We apply Theorem 3.14 and explain in detail in \cite{P-R,DPT},}).
\begin{theorem}\label{Theorem 3.14}\ \\
\begin{enumerate}
\item
[(a)] Let $F=F_{n+1}=\{x_1,x_2,...x_{n+1}:\ \}$ and let 
$F^{(k)}$ be the kernel of the map $F\to Z_k$ which sends $x_i$ to $1$.
Furthermore let $F^{(\infty)}=ker(F\to Z)$.
Define $\bar{F}^{(k)}=F^{(k)}/(x_i^k)$ and  $y_i=x_ix_{n+1}^{-1}$. 
Let $\tau: F \to F$ be an automorphism given by $\tau(w) =x_{n+1}wx_{n+1}^{-1}$.
% and  $y_i=x_ix_{n+1}^{-1}$ and $\bar y_i=xx_ix^{-2}$, 
%where $x\in {x_1,x_2,...,x_n}$ and $x_i\neq x$.
%Then $\{y_i, \bar y_i\}$ forms a basis of $\bar{F}^{(3)}$.
\begin{enumerate}
\item
[(i)] $F^{(k)}$ is a free group generated freely by $nk+1$ elements
$\tau ^j(y_i)$, for $i\leq n$ and $0\leq j \leq k-1$, and $x_{n+1}^k$.
\item [(ii)] 
$F^{(\infty)}$ is freely generated by elements $\tau ^j(y_i)$, for $i\leq n$ 
and any integer $j$.
\item[(iii)] $\bar{F}^{(k)}$ is freely generated by $n(k-1)$ elements 
$\tau ^j(y_i)$, for $i\leq n$ and $0\leq j < k-1$.\\
Notice that one has relations $y_i\tau(y_i)...\tau^{k-1}(y_i)=1$, for any $i$.
\end{enumerate}
\item
[(b)] 
\begin{enumerate} 
\item 
[(i)] ${\Pi}_D^{(k)}\ast \underbrace{Z\ast Z\ast ...\ast Z}_{k-1\ \ times}
= {\Pi}_{D\sqcup O}^{(k)}$ 
has the following Wirtinger type description:\\
There are $k-1$ generators, $\tau ^j(y_i)$, $0\leq j < k-1$, 
corresponding to the $i$'th arc of the diagram $D$, and there are $k-1$
relations $\tau ^j(r_s))$, $0\leq j < k-1$, corresponding to any 
crossing, $v_s$, where $r_s$ depends on the sign of a crossing as follows:
\begin{enumerate}
\item [(+)] In the case of the positive crossing (Fig.3.1(a)):
$r_s=y_i\tau (y_k) (\tau (y_i))^{-1}y_j^{-1}$.
\item [(-)] 
In the case of the negative crossing (Fig.3.1(b)):
$r_s=y_i\tau (y_j) (\tau (y_i))^{-1}y_k^{-1}$.
\end{enumerate}
We have to remember that $\tau^{k-1}(y_i)=
(y_i\tau(y_i)...\tau^{k-2}(y_i))^{-1}$.
%Hint. Generally we need relations of the form $wr_sw^{-1}$ but they can be 
%reduced to $\tau ^j(r_s))$ because $x_ir_sx_i^{-1}= 
%x_ix_{n+1}^{-1}x_{n+1}r_sx_{n+1}^{-1}x_{n+1}x_i^{-1}= y_i\tau(r_s)y_i^{-1}$.
\item 
[(ii)] ${\Pi}_D^{(\infty)}\ast F_{\infty}
=\{\tau^j(y_i),\ i\leq n:\ \tau^j(r_s), j\in Z \}$, where $F_{\infty}$ is a
countably generated free group.
\end{enumerate} 
\end{enumerate}
\end{theorem}

\begin{corollary}\label{Corollary 3.15}
\begin{enumerate}
\item [(i)]  $H_1(M^{(\infty)})\oplus Z[t^{\pm 1}] = Z[t^{\pm 1}]
(y_1,y_2,...,y_n)/(\bar r_s)$, where 
$\bar r_s$ are relations associated to crossings: 
$ (1-t^{\epsilon})y_i+t^{\epsilon}y_k - y_j=0$, where $\epsilon =\pm 1$ is 
the sign of the crossing $v_s$; Fig.3.1.\footnote {We can think of it as
Wirtinger type description of the Burau representation. 
Compare also [Re-1,Ch.II(14)].} 
\item [(ii)]
$H_1(M^{(k)},Z)\oplus (Z)^{k-1}= 
Z[t^{\pm 1}]/(1+t+...t^{k-1})(y_1,y_2,...,y_n)/(\bar r_s)$ 
\item[(iii)] $H_1(M^{(k)},Z_m)\oplus (Z_m)^{k-1}$ has the following ``coloring"
description:\\
Every arc of the diagram is colored by a sequence, 
$(a_0,a_1,...,a_{k-2})$, of colors taken
from the set of $m$ colors, $(0,1,2,...,m-1)$. These colorings form a $Z_m$
module $Z_m^{\lambda (k-1)}$, where ${\lambda}$ is the number of arcs in 
the diagram.
Coloring of an arc can be coded by a polynomial of degree $k-2$ with
coefficients in $Z_m$, $w=\Sigma_{i=0}^{k-2}a_it^i$. Now we consider the space
(submodule) of allowed coloring, that is colorings which at any crossing
satisfy the equation: $ (1-t^{\epsilon})w_i+t^{\epsilon}w_k - w_j=0$, where
$1+t+...t^{k-1}=0$ (in particular $t^{-1}=-1 -t-t^2-...-t^{k-2}$), 
$\epsilon =\pm 1$ as in (i), $w_l$ are polynomials of degree $k-2$ with 
coefficients in $Z_m$ corresponding to arcs at a crossing as in Fig.3.1.
Allowed colorings form a group $H_1(M^{(k)},Z_m)\oplus (Z_m)^{k-1}$,
Compare [S-W].
%Hint. We use the fact that $Hom(P \to Z_m) = P \otimes Z_m$ for a finitely
%generated abelian group $P$.
\end{enumerate}

\end{corollary}
We can generalize the group $H_D$ in yet another direction. We can consider
the $|a|$ by $|v|$  matrix $(b_{i,j})$ where $|a|$ is the number of arcs of 
the diagram and $|v|$ the number of crossings of the diagram, and where 
 $b_{i,j}=0$ if the $i$'th arc is disjoint from the $j$'th crossing,
$b_{i,j}=2$ if the $i$'th arc is 
the overcrossing of the $j$'th crossing and $b_{i,j}=-1$ if the $i$'th arc is
an undercrossing of the $j$'th crossing. We consider the matrix up to
changes caused by Reidemeister moves. Virtually nothing is known about
the invariant given by this matrix (except the group $H_D$). 
One should, at least compare this matrix with the Goeritz matrix and 
the Seifert matrix (see \cite{Gor}).

\section{Coloring and statistical mechanics}\label{4} 

Further modification of our method leads to state models of 
statistical mechanics and the Yang-Baxter equation. 
We will illustrate it by two examples. In the first we consider 
%the Yang-Baxter operator corresponding to n-colorings, 
the state sum corresponding to n-colorings and
in the second 
%we interpret the Yang-Baxter operator leading to the generalized
%Jones polynomial, in terms of weighted colorings.
we give (after Jones) the state sum approach to the skein (generalized Jones) 
polynomial. 
\begin{example}\Label{4.1}
For $n$ colors, $0,1,...,n-1$, every coloring of arcs of a diagram by these
colors, is called a state of the diagram. We associate to any state, $s$, and
any crossing $v$, a weight, $w(v,s)$, depending on the colors of arcs
at $v$. For Fox colorings it will be $1$ if twice the color of the 
overcrossing is congruent to the sum of the colors of the undercrossings 
modulo $n$. It will be $0$ otherwise. We associate to any state, $s$, 
the global weight equal to
${\Pi}_vw(v,s)$. For a Fox coloring it is $1$ if $s$ is an $n$-coloring and
$0$ otherwise. Finally we define the partition function, $Z_D(n)$, 
to be the sum over all states of their global weights, i.e.
$$Z_D(n)={\Sigma}_s{\Pi}_vw(v,s).$$ In our example we get
$Z_D(n)=col_n(D)$. This state sum description of Fox n-colorings
was given in [H-J], and it is called, in statistical mechanics, a vertex type
model. 
\end{example}
In our example, colors were associated to arcs of the diagram; in the general 
case of the vertex model weights are associated to edges of a graph. 
So we have to think
of the link diagram as a graph with additional structure (under-over crossing)
at the vertices. 
\begin{example}\label{4.2}
Let a link diagram be considered as a graph with vertices of valency 4
at crossings and vertices of valency 2 at maxima and minima.
We assume that there is only finite number of extrema  and that crossings
are positioned vertically, so that after smoothing them, one do not
 introduce new extrema.  A state is a
function from edges of the diagram to the set of $k$ colors. We consider
weights to be in $Z[q^{\pm 1}]$. For any vertex, $v$, of degree 4,
and inputs with colors $i,j$ and output with colors $k,l$,
we associate the weight
\[ w_+(i,j;k,l)=\left\{ \begin{array}{lll}
            q-q^{-1} & if & i<j,\ i=k,\ j=l\\
             1       & if & i=l, j=k, i\neq j \\
             q       & if & i=j=k=l  \\
             0       &  otherwise & ,
\end{array}
\right. \] in the case of a positive crossing, and
\[ w_-(i,j;k,l)=\left\{ \begin{array}{lll}
            q^{-1}-q & if & i>j,\ i=k,\ j=l\\
             1       & if & i=l, j=k, i\neq j \\
             q^{-1}  & if & i=j=k=l  \\
             0       &  otherwise & ,
\end{array}
\right. \] in the case of a negative crossing; see Fig.4.1.\\
For any vertex of degree 2, we associate the weight 
$w(i)=q^{\pm 1/2(2i-k-1)}$, according to the convention of Fig.4.2.\\
As proven by Jones [Jo-1] the partition function given by this vertex 
model is equal to the version of the skein (homflypt) polynomial for 
regular isotopy classes of links (associated with the skein relation 
$q^kP_{L_+} -q^{-k}P_{L_-}=(q-q^{-1})P_{L_0}$).
\end{example}

%\vspace*{2in}
\centerline{\psfig{figure=states.eps}}
\begin{center}
Fig. 4.1 
\end{center}

\ \\ 
%\vspace*{2in}
\centerline{\psfig{figure=max-min.eps}}
\begin{center}
Fig. 4.2
\end{center}

%{\bf Bibliography}
%\ \\

\ \\
 \ \\

\ \\
Department of Mathematics\\
University of California\\
Berkeley, CA 94 720\\
e-mail: Jozef@math.berkeley.edu

Current address:\\
Department of Mathematics\\
George Washington University  \\
Washington, DC 20052 \\
e-mail: przytyck@gwu.edu
\end{document}